\documentclass[12pt,a4paper]{article}

\usepackage{latexsym,amssymb,amsthm}
\usepackage{url,graphics}
\usepackage[dvips]{graphicx,epsfig}

\title{\bf  An Upper bound on the chromatic number of circle graphs without $K_4$}
\author{
     \begin{tabular}{c}
     G.\,V.\,Nenashev
     \\[6pt]
     \\[-2pt]
     {\small E-mail: \texttt{glebnen@mail.ru}       }
    \end{tabular}
    }
\date{}

\begin{document}
\maketitle
\righthyphenmin=2
\renewcommand*{\proofname}{\bf Proof}
\newtheorem{thm}{Theorem}
\newtheorem{lem}{Lemma}
\newtheorem{cor}{Corollary}
\theoremstyle{definition}
\newtheorem{defin}{Definition}
\theoremstyle{remark}
\newtheorem*{rem}{\bf Remark}

\def\I{{\rm Int}}
\def\R{{\rm Bound}}
\def\q#1.{{\bf #1.}}
\renewcommand\geq{\geqslant}
\renewcommand\leq{\leqslant}

\abstract{Let $G$ be a circle graph without clique on 4 vertices.
We prove that the chromatic number of  $G$  doesn't exceed 30.
}

\section{\bf Introduction}

\begin{defin}

Fix a circle. Let $\cal A$ be a finite set of chords of this circle. 
A {\it circle graph} $G(\cal A)$ is a graph with $\cal A$ as a vertex set, in 
which two vertices
are adjacent if and only if corresponding chords have a common inner point.
\end{defin}

\begin{figure}[htb!]
\centering
\includegraphics[scale=0.25]{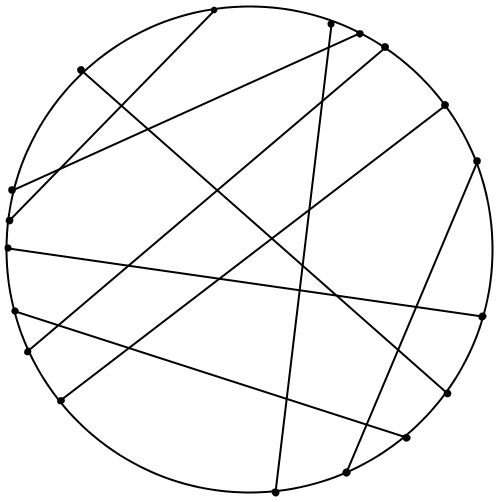}
\caption{A circle graph}
\end{figure}

As usual, let a {\it proper coloring} of a graph be such a coloring of vertices 
of this graph, that any two adjacent vertices have different
colors.

We denote by~$\omega(G)$ the size of a maximal clicque of the graph~$G$ 
and by~$\chi(G)$ the {\it chromatic number} of the graph~$G$,
i.e. the minimal number of colors in a proper coloring of this graph.

Let~$G$ be a circle graph. One can find several results about bounds of~$\chi(G)$.
In~1988 A.~V.~Kostochka~\cite{K1}  proved that $\chi(G)\leq{5}$ if $\omega(G)=2$.
In~1996 A.~A.~Ageev~\cite{A96}  introduced an example of a circle graph~$G$ with~$\omega(G)=2$ and~$\chi(G)={5}$, i.e. proved that 
the aforementioned bound on the chromatic number  for~$\omega(G)=2$ is tight.
In~1997 A.~V.~Kostochka and J.~Kratochvil~\cite{KK} proved, that $\chi(G)\leq {2^{\omega(G)+6}}$.
In~1999 A.~A.~Ageev~\cite{A99} proved, that if $\omega(G)=2$  and the graph 
$G$ has girth at least~5 (i.e. contains no cycles of length~3 and~4), then~$\chi(G)\leq{3}$.

In~2011 A.~V.~Kostochka and K.~G.~Milans~\cite{KM} proved that a circle graph without~$K_4$ has a proper coloring in~38 colors.
We prove that such a graph has a proper coloring in~30 colors.

\section{\bf Proper coloring  of a circle graph without~$K_4$ in 30 colors.}

Hereafter we mean by a  {\it coloring of  chords of the set~$\cal A$}  
a proper coloring of the vertices of the graph~$G(\cal A)$.
We work with arcs and chords of a fixed circle.

\begin{lem}
\label{lf3}
Let an arc $(X,Y)$ and two sets of chords~$\cal A$ and~$\cal B$  satisfy the following conditions:

{\bf (1)} Each triangle in  $G(\cal A \cup \cal B)$ contains at most one chord of the set~$\cal B$.

{\bf (2)} Each chord  of the set $\cal A$ has exactly one end on the arc~$(X,Y)$.

{\bf (3)} For each chord~$b_j \in \cal B$ both ends of $b_j$ are on the arc~$(X,Y)$, and there exists a chord~$a_i\in \cal A$, such
that~$a_i$ intersect~$b_j$.

Then there is a proper coloring of the set~of chords $\cal B$ in~$3$ colors.

\begin{figure}[htb!]
\label{f2}
\centering
\includegraphics[scale=0.25]{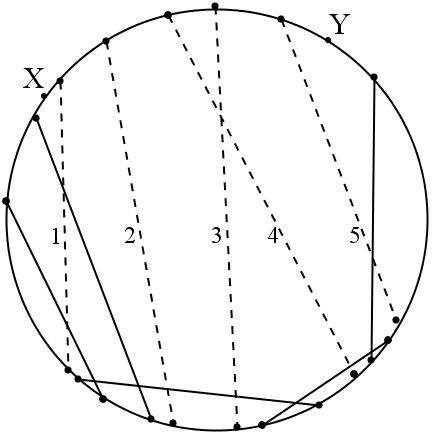}
\caption{Sets of chords~$\cal A$ and~$\cal B$}
\end{figure}

\begin{proof}

We suppose that chords of the set~$\cal A$ are enumerated in the order in 
which their ends lie on the arc~$(X,Y)$ (from~$X$ to~$Y$).
On the figure~2 chords of  the set~$\cal A$ are  shown as  
dashed lines and are enumerated, chords of the set~$\cal B$ are  shown as solid lines.

For each chord~$b_j\in \cal B$ denote by~${P_j+1}$ and~$Q_j$ the minimal and maximal
numbers  of chords of~$\cal A$, that intersect~$b_j$. Consider the interval~$P_jQ_j$ on the real
line.

Note that our intervals cannot intersect each other, i.e. the situation when $P_i<P_j<Q_i<Q_j$ for two  intervals~$P_iQ_i$ and~$P_jQ_j$
is impossible. Otherwise correspondent chords and the chord~$a_{P_j}$ form a triangle in the graph~$G(\cal A \cup \cal B)$.

If two chords from~$\cal B$ intersect each other, then corresponding intervals 
will {\it touch} each other, i.e. the left end of one of them coincides with the right end of another.
If we obtain equal intervals for some two chords, we identify them and consider such interval only once.

Now we shall color these intervals in 3 colors in such a way that any two touching intervals have different colors.
It is easy to see that the corresponding coloring of chords is proper.

We shall color each interval and its left end in the same color (i.e. the color of an interval is the color of its left end).
Let us describe the order of coloring. On each step we do the following.

$1^\circ$ If it is possible, we choose the   {\it maximal} uncolored left end, which is a right end of a colored interval.

$2^\circ$ If we cannot perform item $1^\circ$, then we choose the {\it minimal}  left end among uncolored ends.

Let's prove that on each step we can find a color for the choosen left end.
Let~$L$ be this left end and right ends~$R_1,\dots,R_n$ form 
intervals with~$L$ (right ends are enumerated in increasing order).

Let's prove that at most one of the points~$R_1,...,R_n$ is already colored as a left end.  Assume the contrary, let~$R_i$ and~$R_j$ 
be colored, where~$i<j$.  Consider the first colored point~$D$ strictly inside the interval~$LR_j$.
Consider the moment when we have colored~$D$. At that moment there was no colored interval~$ID$ with right end~$D$, since otherwise
either  $I$ lies on~$LR_j$ and was colored before~$D$, or~$ID$ intersects~$LR_j$.  Thus we cannot color~$D$ according to item~$1^\circ$.
But we also cannot color~$D$ according to item~$2^\circ$, since $L<D$ and~$L$ was not colored in the considered moment.
We obtain a contradiction.

Let's prove that at most one interval with end~$L$ is colored.
Suppose the contrary, let two intervals~$P_1L$ and~$P_2L$ 
be colored (since $L$ is not colored, then~$L$ is the right end 
of both colored intervals). Assume that~$P_1$ was colored before~$P_2$.
At the moment of coloring~$P_2$ the point~$P_1$ was colored (and, hence, the interval $P_1L$ was colored).
 Since the point~$L$ is not colored and is the right end of a colored interval, we must color accordingly
to item~$1^\circ$. But in this case we would color~$L$ instead of~$P_2$,
since~$L > P_2$. We obtain a contradiction.

Thus not more than two colors are forbidden for~$L$ and we can color~$L$.

Hence all chords of the set~$\cal B$ can be proper colored in  3 colors according to 
the coloring of correspondent intervals.
\end {proof}

\end{lem}

\begin{lem}
\label{t1f4}
Let an arc $(X,Y)$ and two sets of chords~$\cal A$ and~$\cal B$  satisfy the following conditions:

 {\bf (1)} The graph $G(\cal A \cup \cal B)$ does not contain~$K_4$.

{\bf (2)} Each chord  of the set $\cal A$ has at most one end on the arc~$(X,Y)$.

{\bf (3)} For each chord~$b_j \in \cal B$ both ends of $b_j$ are on the arc~$(X,Y)$, and there exists such a chord~$a_i\in \cal A$,
that~$a_i$ intersect~$b_j$.

Then there is a proper coloring of chords of the set~$\cal B$ in~$15$ colors.

\begin{proof}
Note, that chords of the set~$\cal A$, which have no end on~$(X,Y)$, cannot intersect 
chords of the set~$\cal B$. We remove all such chords from~$\cal A$. Now 
every chord of~$A$ has exactly one end on~$(X,Y)$.

We suppose that chords of the set~$\cal A$ are enumerated in the same order
as their ends on the arc~$(X,Y)$ (from~$X$ to~$Y$).
Let $A_j$ be the end of a chord~$a_j \in \cal A$ on the arc~$(X,Y)$. 
We add points~$A_0=X$, $A_{last}=Y$.

We construct a new set of chords~$\cal C$. For each~$j$ consider 
all ends of chords from~$\cal B$, that lie between~$A_j$ and~$A_{j+1}$.
We change the order of these ends. 
At first we mark all the right ends (with the same enumeration as before),
and then mark all the left ends (also in the same order as before). 
Since no chord has both ends between $A_j$ and $A_{j+1}$, such transformation is possible.
An example of this transformation is shown on figure~3.

{\begin{figure}[htb!]
\centering
\includegraphics[scale=0.3]{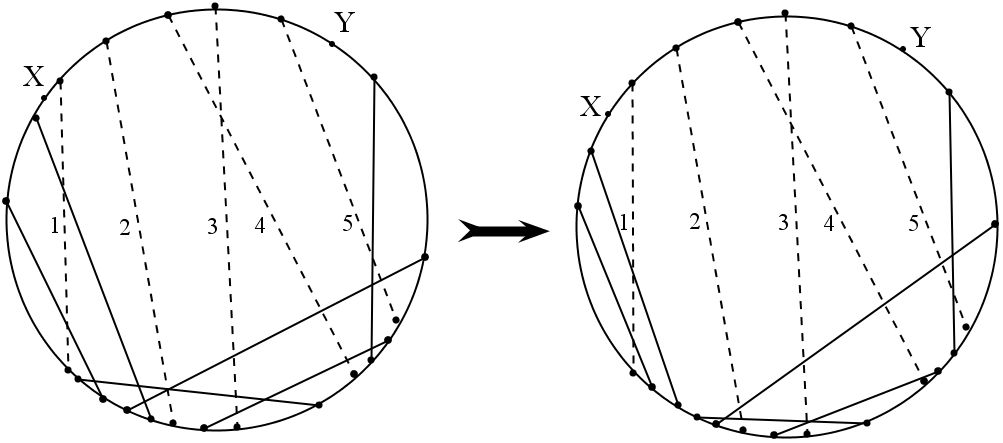}
\caption{Construction of the set $\cal C$}
\end{figure}
}

In the constructed set  there are no intersections of pairs of chords, one of which 
have  the right end between  $A_j$ and~$A_{j+1}$ and the other
have the  left end between  $A_j$ and~$A_{j+1}$.
All other intersections of the chords of the set $ \cal A \cup \cal B $ remain in the set $ \cal A \cup \cal C $.

Let's prove that there is no triangle in  $G({\cal C})$. 
Assume the contrary, let chords $L_1R_1$, $L_2R_2$ and~$L_3R_3$ form a triangle.
Then by construction of the set $\cal C$ there is a point~$A_j$ between 
$\max(L_1,L_2,L_3)$ and $\min(R_1,R_2,R_3)$, hence there 
is a subgraph~$K_4$ in $G({\cal A} \cup {\cal C} )$. 
Then  there is a subgraph~$K_4$ in $G({\cal A} \cup {\cal B})$, we get a contradiction.

Thus, there is no triangle in  $G({\cal C})$ and then by the work~\cite{K1}  
we can color chords of the set~$\cal C$ in  5 colors. We consider the 
coloring of chords of the set $\cal B $ in 5 colors, corresponding to a 
proper  coloring  of chords  of~$\cal C$.

Let us correct the coloring of chords of~$\cal B$. Consider a color~$w$ and the set of chords~${\cal B}_w$ which have this color.
Let's prove that the graph~$G({\cal B}_w \cup \cal A)$  cannot have a triangle with 
at least 2 chords of~${\cal B}_w$.
By the construction of the set~$\cal C$, if two chords~$b_k,b_\ell\in {\cal B}_w$ 
intersect each other, then the left end of one of them and 
the right end of the other lie between $A_i$ and $A_{i+1}$ for some~$i$, see figure~4.

{\begin{figure}[htb!]
\centering
\includegraphics[scale=0.25]{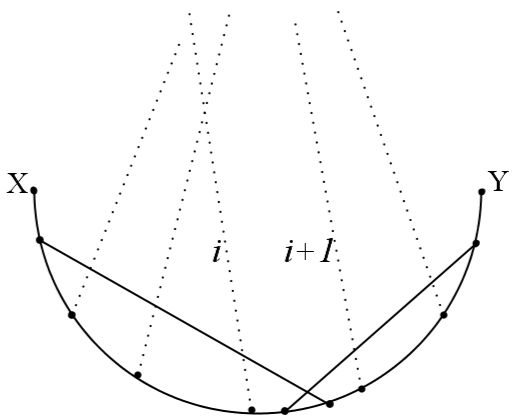}
\caption{Intersecting chords of~${\cal B}_w$ }
\end{figure}
}

Hence no chord of the set ~$\cal A$   can intersect both chords 
$b_k$ and~$b_\ell$, consequently, there is no triangle with 
two chords of~${\cal B}_w$ and one chord of~$\cal A$. If some chord 
of~${\cal B}_w$ intersects both chords~$b_k$ and~$b_\ell$, then 
this chord intersects  $a_i$ or $a_{i+1}$,  thus there is a triangle with 
two chords of~${\cal B}_w$ and one chord of~$\cal A$. But we have already proved
that it is not possible.

Then by lemma~\ref{lf3} we can color chords of the set~ ${\cal B}_w$ in 3 colors. 
Thus we can color all chords of the set~${\cal B}$  in $15= 5\cdot 3$ colors and this coloring will be proper.
\end{proof}
\end{lem}

\begin{lem}
\label{t2f4}
Let an arc $(X,Y)$ and two sets of chords~$\cal A$ and~$\cal B$  satisfy the following conditions:

{\bf (1)} The graph $G(\cal A \cup \cal B)$ does not contain~$K_4$.

{\bf (2)} Each chord  of the set $\cal A$ has at most one end on the arc~$(X,Y)$.

{\bf (3)} For each chord~$b_j \in \cal B$ both ends of $b_j$ are on the arc~$(X,Y)$.

{\bf (4)} Every connected component of the graph  $G(\cal A \cup \cal B)$ contains 
at least one chord of the set~$\cal A$.

Assume additionally that there exists a proper coloring of the set~$\cal A$ in (at most)~$15$ colors.
Then it is possible to color chords of the set~$\cal B$ such that the resulting coloring of the set~${\cal A}\cup {\cal B}$
is a proper coloring in not more than~$30$ colors.

\begin{proof}

Note that chords of the set~$\cal A$, which have no end on~$(X,Y)$, cannot intersect 
chords of the set~$\cal B$. We remove all such 
chords from~$\cal A$. Now any chord of~$A$ has exactly one end on~$(X,Y)$.

We suppose that chords of the set~$\cal A$ are enumerated in the order in which 
their ends lie on the arc~$(X,Y)$ (from~$X$ to~$Y$).
Let $A_j$ be the end of a chord~$a_j \in \cal A$ on the arc~$(X,Y)$. We add points~$A_0=X$, $A_{last}=Y$.

We prove the statement of Lemma by induction on $|\cal B|$.  The base is the obvious case $|\cal B|$=0.

{\sf Induction step.}
Let $\cal C$ be the set of all chords of~$\cal B$, that 
intersect chords of~$\cal A$ and~${\cal B}' = {\cal B} \setminus {\cal C}$. Clearly,  
$\cal C$ is nonempty. 
By Lemma~\ref{t1f4} there is a proper coloring of chords of the set~$\cal C$ in 15 colors.
We use for this purpose 15 colors not used 
for coloring of the set $\cal A$ and obtain a proper 
coloring of the graph $G({\cal A} \cup \cal C)$ in not more than 30 colors.

Since no chord of $\cal B'$ intersect a chord of~$\cal A$,  then for each 
chord  $b \in \cal B'$ there exists such~$j$, that both ends 
of the chord~$b$ lie between~$A_j$ and~$A_{j+1}$. Let~${\cal B}_j$ 
be a set of all chords from~$\cal B'$ such that both ends of these chords 
lie between~$A_j$ and~$A_{j+1}$. Then ${\cal B}'=  \bigcup {\cal B}_j$.
Obviously,  chords of different sets~${\cal B}_j$ do not intersect each other.

It is enough to color each set of chords ${\cal B}_j$ separately so that the coloring of the graph
$G({\cal A} \cup {\cal C} \cup {\cal B}_j)$ is proper. 
Since no chord of ${\cal B}_j$ intersect a chord of~$\cal A$ and the coloring of
 $G({\cal A} \cup \cal C)$ is proper, it suffices to color ${\cal B}_j$ 
 so that the coloring of the graph
$G({\cal C} \cup {\cal B}_j)$ is proper.

Let us consider an arc~$(A_j,A_{j+1})$ and prove by induction, 
that we can color chords of the set~${\cal B}_j$. We
 use the arc~$(A_j,A_{j+1})$ instead of the arc~$(X,Y)$, and sets of 
 chords $\cal C$ and~${\cal B}_j$ instead of~$\cal A$ and~$\cal B$ . Let us verify the conditions of Lemma.

1) The graph $G({\cal C} \cup {\cal B}_j)$ is a subgraph of the graph~$G(\cal A \cup \cal B)$, hence it does not contain~$K_4$.

2)  Each chord  of the set $\cal C$ has at most one end on the arc~$(A_j,A_{j+1})$, since this chord intersects some chord from~$\cal A$.

3) This condition is clear by the construction of~${\cal B}_j$.

4) Assume that the condition~$(4)$ does not hold and there is a connected component in the graph $G({\cal C} \cup {\cal B}_j)$,
which contais no chords of~$\cal C$. Chords of this component do not intersect chords of the set~$\cal A$
(since only chords of~$\cal C$ intersect chords of~$\cal A$) and do not intersect chords of any set ${\cal B}_i$ for $i\ne j$.
Hence this connected component (without chords of the set $\cal A$) is a connected component of the graph $G(\cal A \cup \cal B)$,
that contradicts to the condition of Lemma.

Thus  conditions hold and $|{\cal B}_i|<|\cal B|$ (otherwise $|{\cal C}|=0$). 
Now we can use the induction assumption and finish the proof.
\end{proof}
\end{lem}

\begin{thm}
\label{t3f4}
Let the set of chords $\cal A$ be such that the graph $G(\cal A)$ does not contain~$K_4$.
Then there is a proper coloring of chords of the set~$\cal A$ in~$30$ colors.

\begin{proof}
We assume, that the circle graph~$G(\cal A)$ is connected (otherwise we color each connected component independently).

Consider a chord $CD$ from $\cal A$ and any arc $(X,Y)$, which contains all ends of the chords of the set~$\cal A$, except~$D$.
Let ${\cal B}={\cal A} \setminus  \{CD\}$. We color the chord~$CD$ in color~1.

Let us verify the conditions of lemma~\ref{t2f4} for the arc~$(X,Y)$, first set of chord~$\{CD\}$ and second set of chords~$\cal B$.

1) The graph~$G(\{CD\}\cup {\cal B})=G(\cal A)$  does not contain~$K_4$.

2) Exactly one end of the chord~$CD$ lies on the arc~$(X,Y)$.

3) Both ends of any chord from~$\cal B$ lie on the arc~$(X,Y)$.

4) Since the graph~$G(\{CD\} \cup {\cal B})=G(\cal A)$ is connected, it has only one 
connected component and this component contains~$CD$.

Thus all conditions hold and we can color chords of the set $\cal A$ in  30 colors.
\end{proof}
\end{thm}

\begin{rem}
By the same method one may show that any circle graph without~$K_n$ has a proper 
coloring in~$5\cdot 6^{n-3} $ colors.
\end{rem}

Translated by D.\,V.\,Karpov, edited by F.\,V.\,Petrov.


\begin{thebibliography}{9}



\bibitem{K1} {\sc  A.\,V.\,Kostochka.} \textit{On upper bounds for the chromatic numbers of graphs.}
Trudy Instituta Mathematiki {\bf 10}, p.\,204-226, 1988.


\bibitem{A96} {\sc  A.\,A.\,Ageev.} \textit{A triangle-free circle graph with chromatic number $5$.}
Discrete Math. {\bf 152}, p.\,295-298, 1996.



\bibitem{KK} {\sc  A.\,V.\,Kostochka, J.\,Kratochvil.} \textit{Covering and coloring polygon-circle graphs.}
Discrete Math. {\bf 163}, p.\,299-305, 1997.

\bibitem{A99} {\sc  A.\,A.\,Ageev.} \textit{Every circle graph of girth at least $5$ is $3$-colourable.}
Discrete Math. {\bf 195}, p.\,229-233, 1999.

\bibitem{K2} {\sc  A.\,V.\,Kostochka.}
\textit{Coloring intersection graphs of geometric figures with a given clique number.}
Contemp. Math. {\bf 342}, p.\,127-138, 2004.

\bibitem{JC} {\sc J.\,Cerny.} \textit{Coloring circle graphs.}
Electronic Notes in Discrete Mathematics {\bf 29}, p.\,457-461, 2007.

\bibitem{KM} {\sc  A.\,V.\,Kostochka, K.\,G.\,Milans}
\textit{Coloring clean and $K_4$-free circle graphs.}  2011, submitted.



\end{thebibliography}
\end{document}